\documentclass[12pt,reqno]{amsart}

\usepackage{amssymb}

\usepackage{eucal}

\input{diagrams}

\font\sss=cmss8

\def\cL{{\mathcal L}}
\def\cM{{\mathcal M}}
\def\cN{{\mathcal N}}
\def\cO{{\mathcal O}}

\def\BN{{\mathbb N}}

\def\BP{{\mathbb P}}

\def\BZ{{\mathbb Z}}

\def\sD{\mbox{\sf D}}

\def\sK{\mbox{\sf K}}

\def\b{\operatorname{b}}

\def\Chain{\mbox{\sf Ch}}

\def\Coh{\mbox{\sf coh}}

\def\colim{\operatorname{colim}}
\def\D{\sD}
\def\Db{\sD^{\operatorname{b}}}

\def\dim{\operatorname{dim}}

\def\Ext{\operatorname{Ext}}

\def\gldim{\operatorname{gldim}}
\def\gr{\mbox{\sf gr}}
\def\Gr{\mbox{\sf Gr}}

\def\h{\operatorname{h}}

\def\Hom{\operatorname{Hom}}

\def\opp{\operatorname{op}}

\def\Proj{\mbox{\sf Proj}}
\def\QCoh{\mbox{\sf QCoh}}

\def\RHom{\operatorname{RHom}}

\def\tails{\mbox{\sf qgr}}
\def\Tails{\mbox{\sf QGr}}

\def\Tors{\mbox{\sf Tors}}
\def\Tot{\operatorname{Tot}}

\def\Z{\operatorname{Z}}

\numberwithin{equation}{part}

\newtheorem{Lemma}{Lemma}[section]
\newtheorem{Theorem}[Lemma]{Theorem}
\newtheorem{Proposition}[Lemma]{Proposition}

\theoremstyle{definition}
\newtheorem{Definition}[Lemma]{Definition}
\newtheorem{Setup}[Lemma]{Setup}

\newtheorem{Remark}[Lemma]{Remark}

\newtheorem{Example}[Lemma]{Example}

\def\Nc{Non-com\-mu\-ta\-ti\-ve}
\def\nc{non-com\-mu\-ta\-ti\-ve}
\def\BGGc{BGG correspondence}
\def\ASreg{AS regular}

\def\A{A}                
\def\Ad{A^!}             
\def\Ao{\A^{\opp}}

\def\Ealg{E}             

\def\calO{\cO}           

\def\lm{left-mo\-du\-le}

\def\Abm{$\A$-bi-mo\-du\-le}
\def\Alm{$\A$-left-mo\-du\-le}

\def\Adlm{$\Ad$-left-mo\-du\-le}

\def\DGrA{\D(\Gr\,\A)}

\def\DTailsA{\D(\Tails\,\A)}
\def\DbTailsA{\D^{\b}(\Tails\,\A)}

\def\GrInj{\mbox{\sf Inj}}
\def\GrTFInj{\mbox{{\sf Inj}$\,{}_{\operatorname{tf}}$}}
\def\TailsInj{\mbox{\sf QInj}}
\def\GrFree{\mbox{\sf Free}}
\def\GrCoFree{\mbox{\sf CoFree}}
\def\GrStab{\overline{\mbox{\sf Gr}}}
\def\grstab{\overline{\mbox{\sf gr}}}

\def\scrDGrA{\mbox{\sss D}(\mbox{\sss Gr}\,\A)}
\def\scrDTailsA{\mbox{\sss D}(\mbox{\sss QGr}\,\A)}
\def\scrGr{\mbox{\sss Gr}}
\def\scrGrCoFree{\mbox{\sss CoFree}}
\def\scrGrStab{\overline{\mbox{\sss Gr}}}
\def\scrsK{\mbox{\sss K}}
\def\scrQCoh{\mbox{\sss QCoh}}
\def\scrTails{\mbox{\sss QGr}}

\def\Romega{\operatorname{R}\!\omega}
\def\BGGF{\operatorname{F}}
\def\BGGG{\operatorname{G}}
\def\BGGphi{\operatorname{\varphi}}
\def\BGGgamma{\operatorname{\gamma}}

\def\elliptic{C}	
\def\point{p}		
\def\pointtwo{q}	
\def\pointmodule{P}	
\def\auto{\tau}		
\def\bundle{\cL}	

\begin{document}

\title[Non-commutative BGG]
{A non-commutative BGG correspondence}

\author{Peter J\o rgensen}
\address{Danish National Library of Science and Medicine, N\o rre
All\'e 49, 2200 K\o \-ben\-havn N, DK--Denmark}
\email{pej@dnlb.dk, www.geocities.com/popjoerg}




\begin{abstract} 

A \nc\ version of the
Bern\v{s}te\u{\i}n-Gel${}^\prime$fand-Gel${}^\prime$fand (BGG)
correspondence is set up, and a sample application is given to
periodic injective resolutions.

\end{abstract}

\maketitle

\setcounter{section}{-1}
\section{Introduction}
\label{sec:introduction}

The Bern\v{s}te\u{\i}n-Gel${}^\prime$fand-Gel${}^\prime$fand 
correspondence is surprising.  It gives an equivalence of categories
\[
    \grstab(\Ealg) 
    \simeq \Db(\Coh\,\BP^{d-1}).
\]
This was proved in \cite[thm.\ 2]{BGG} and as always, to explain such a
formula requires lots of words.  On the left hand side, $\Ealg$ is the
exterior algebra $\bigwedge(Y_1,\ldots,Y_d)$ and $\grstab(\Ealg)$ is
the category of finitely generated graded $\Ealg$-\lm s modulo
morphisms which factor through injectives.

On the right hand side, $\BP^{d-1}$ is $(d-1)$-dimensional projective
space, $\Coh\,\BP^{d-1}$ is the category of coherent sheaves on
$\BP^{d-1}$, and $\Db(\Coh\,\BP^{d-1})$ is the derived category of
bounded complexes of such sheaves.  The surprising thing about the
Bern\v{s}te\u{\i}n-Gel${}^\prime$fand-Gel${}^\prime$fand (BGG)
correspondence is that the geometric object on the right hand side is
equivalent to the purely algebraic object on the left hand side.

Put differently, if one did not know about the \BGGc, it would really
not be obvious that it is possible to recover $\Db(\Coh\,\BP^{d-1})$
purely algebraically!

In this paper, I will generalize the \BGGc\ to \nc\ projective
geometry.  \Nc\ projective geometry is well established; one of the
seminal papers is \cite{ArtinZhang} but many have been published
since, showing how a range of projective geometry can be generalized
in a \nc\ way.  This turns out also to be true of the \BGGc\ which is
generalized in theorem \ref{thm:BGG} below, and now takes the form
\[
    \GrStab(\Ad) \simeq \DTailsA.
\]
Here $\A$ is a suitable \nc\ graded algebra with Koszul dual algebra
$\Ad$, and the category $\Tails(\A)$ is a \nc\ analogue of the
category $\QCoh(\BP^{d-1})$ of quasi-coherent sheaves on $\BP^{d-1}$. 

After proving this, I consider a sample application to periodic
injective resolutions.  The background is Eisenbud's result
\cite[thm.\ 2.2]{Eisenbud}:  Let $M$ be a finitely generated graded
module without injective direct summands over the exterior algebra
$\Ealg$, for which the Bass numbers $\mu^i(M) = \dim_k
\Ext_{\Ealg}^i(k,M)$ are bounded for $i \geq 0$.  Then the minimal
injective resolution $I$ of $M$ is periodic with period one: All the
$I^i$ and all the $\partial_I^i$ are the same, up to isomorphism and
degree shift.  (In fact, Eisenbud worked with minimal free
resolutions, but using the Matlis duality functor $\Hom_k(-,k)$ on his
result gives the above.)  In section \ref{sec:periodicity}, I start by
showing that this phenomenon can be understood geometrically in a very
simple way:

Using the \BGGc, the module $M$ can be translated to a geometric
object on $\BP^{d-1}$.  Since the Bass numbers of $M$ are bounded,
this object turns out to have zero dimensional support, so is stable
under twisting, that is, tensoring by $\cO_{\BP^{d-1}}(1)$.
Translating back, this means that $M$ is its own first syzygy, and
periodicity of the minimal injective resolution follows.

Next in section \ref{sec:periodicity}, I consider the \nc\ case where
a similar procedure yields remarkably different results: Let $\A$ be a
suitable \nc\ graded algebra, and let $M$ be a finitely generated
graded module over the Koszul dual $\Ad$, for which the Bass numbers
$\mu^i(M)$ are bounded for $i \geq 0$.  Then, choosing $\A$ and $M$
prudently, it is possible to make the minimal injective resolution of
$M$ periodic with any finite period, or to make it aperiodic.

The reason is that when translating $M$ through the \nc\ \BGGc, one
still obtains a geometric object with zero dimensional support.
However, due to the \nc\ (hence non-local) nature of the situation, it
is no longer true that such an object is invariant under twisting.
Rather, the object can have orbit of any finite length, or have
infinite orbit.  Translating back gives the above results on
periodicity of the minimal injective resolution.

Note that the concrete example I will give of this behaviour is
already known from \cite{Smith}.  But the present geometric view
through the \BGGc\ is new.

After these remarks, let me end the introduction with a synopsis of
the paper.  Section \ref{sec:Gr_and_Tails} exhibits $\DTailsA$ as a
full subcategory of $\DGrA$.  Section \ref{sec:Koszul} considers a
version of Koszul duality.  Section \ref{sec:BGG} combines these
results into the \nc\ \BGGc, and shows that under the correspondence,
the simple module $k$ over $\Ad$ corresponds to the ``structure
sheaf'' $\calO$ in $\DTailsA$.

Section \ref{sec:computations} does a few computations which are put
to use in section \ref{sec:periodicity}, where the \BGGc\ is applied
to periodicity of minimal injective resolutions.

To avoid a lengthy section on nomenclature, hints on notation are
given along the way.  The reader should rest assured that no new, let
alone revolutionary, notation is introduced.  The paper remains firmly
on classical ground, and differs notationally only in minor details
from such papers as \cite{ArtinZhang}, \cite{PJCM}, and \cite{Smith}.
However, I do need to perform the following blanket setup which
applies throughout.

\begin{Setup}
\label{set:blanket}
$k$ is a field, and $\A = k \oplus \A_1 \oplus \A_2 \oplus \cdots$ is
a connected $\BN$-graded noetherian $k$-algebra which is \ASreg\ and
Koszul (see \cite[p.\ 206]{PJCM} and \cite[def.\ 1.2.1]{BGS}, or
remark \ref{rmk:blanket}).  I assume $\gldim \A = d \geq 2$.

By \cite[cor.\ 2.3.3]{BGS}, the algebra $\A$ is quadratic, that is,
it has the form 
\[
  \A \cong \mbox{T}(V)/(R)
\] 
where $V$ is a finite dimensional vector space, $\mbox{T}(V)$ the
tensor algebra, and $(R)$ the two sided ideal generated by a space of
relations $R$ in $V \otimes_k V$.  Let $(-)^{\prime}$ denote
$\Hom_k(-,k)$ and define $R^{\perp}$ by the exact sequence
\[
  0 \rightarrow R^{\perp} 
  \longrightarrow V^{\prime} \otimes_k V^{\prime} 
  \longrightarrow R^{\prime} \rightarrow 0.
\]
Then the Koszul dual algebra of $\A$ is
\[
  \Ad = \mbox{T}(V^{\prime})/(R^{\perp}),
\]
see \cite[dfn.\ 2.8.1]{BGS}.
\end{Setup}

\begin{Remark}
\label{rmk:blanket}
\noindent
(i)  For $\A$ to be \ASreg\ means that $\gldim \A = d$ is finite, and
that the graded \Abm\ $k = \A / \A_{\geq 1}$ satisfies
\[
  \Ext_{\A}^i(k,\A) \cong \Ext_{\Ao}^i(k,A) \cong
  \left\{
    \begin{array}{cl}
      0       & \mbox{ for } i \not= d, \\
      k(\ell) & \mbox{ for } i = d
    \end{array}
  \right.
\]
for some $\ell$.  As usual, $(-)(\ell)$ denotes $\ell$'th degree shift
of graded modules, so $M(\ell)_i = M_{i + \ell}$.

\smallskip
\noindent
(ii)  For $\A$ to be Koszul means that the minimal free resolution $L$
of the graded \Alm\ $k = \A / \A_{\geq 1}$ is linear.  That is, the
$i$'th module $L_i$ has all its generators in graded degree $i$, so
has the form $\coprod \A(-i)$.

\smallskip
\noindent
(iii)  It is easy to see that since $\A$ is Koszul, the constant $\ell$
in (i) must be $d$.

\smallskip
\noindent
(iv) By \cite[thm.\ 2.10.1]{BGS} there is the isomorphism
$(\Ad)^{\opp} \cong \Ext_{\A}(k,k)$.  Combining this with $\gldim \A =
d$ gives that $\Ad$ is concentrated in graded degrees $0, \ldots, d$.

\smallskip
\noindent
(v)  The algebra $\Ad$ is graded Frobenius by \cite[prop.\ 5.10]{Smith}.
This means that $\dim_k \Ad$ is finite, and that there is an
isomorphism of graded \Adlm s $(\Ad)^{\prime} \cong \Ad(m)$, where
$(\Ad)^{\prime} = \Hom_k(\Ad,k)$ is the Matlis dual module of $\Ad$.

\smallskip
\noindent
(vi)  Since $\Ad$ is concentrated in graded degrees $0, \ldots, d$, the
constant $m$ in (v) must be $d$.  So there is an isomorphism of graded
\Adlm s $(\Ad)^{\prime} \cong \Ad(d)$.
\end{Remark}

\section{The categories $\Gr(\A)$ and $\Tails(\A)$}
\label{sec:Gr_and_Tails}

\begin{Remark}
\label{rmk:Romega}
Let me first recapitulate a few items from \cite{ArtinZhang}, to which
I refer for further details and proofs.

The category $\Gr(\A)$ has as objects all $\BZ$-graded \Alm s and as
morphisms all homomorphisms of \Alm s which preserve graded degree.

A module $M$ in $\Gr(\A)$ is called torsion if each $m$ in $M$ is
annihilated by $\A_{\geq n}$ for some $n$.  The torsion modules form a
dense subcategory $\Tors(\A)$ of $\Gr(\A)$, and the quotient category
is
\[
  \Tails(\A) = \Gr(\A) / \Tors(\A).
\] 
This category behaves like the category of quasi-coherent sheaves on
the space $\Proj(\A)$, although $\Proj(\A)$ itself may not make sense.
For instance, if $\A$ is commutative, then $\Tails(\A)$ is in fact
equivalent to the category of quasi-coherent sheaves on $\Proj(\A)$ by
Serre's theorem, as given in \cite[thm., p.\ 229]{ArtinZhang}.

The degree shifting functor $(-)(1)$ on $\Gr(\A)$ induces a
functor on $\Tails(\A)$ which I will also denote $(-)(1)$.

The category $\Gr(\A)$ has the full subcategory $\gr(\A)$ consisting
of finitely generated modules.  Induced by this, $\Tails(\A)$
has the full subcategory $\tails(\A)$ which behaves like the category
of coherent sheaves on $\Proj(\A)$.

The projection functor $\Gr(\A) \stackrel{\pi}{\longrightarrow}
\Tails(\A)$ has a right-adjoint functor $\Tails(\A)
\stackrel{\omega}{\longrightarrow} \Gr(\A)$ by \cite[p.\
234]{ArtinZhang}, so there is an adjoint pair
\[
  \begin{diagram}[labelstyle=\scriptstyle,midshaft,width=1.75cm]
    \Gr(\A) & \pile{ \rTo^{\pi} \\
                     \lTo_{\omega} } & \Tails(\A). \\
  \end{diagram}
\]
As follows from \cite[prop.\ 7.1]{ArtinZhang}, these functors send
injective objects to injective objects, and restrict to a pair of
quasi-inverse equivalences
\begin{equation}
\label{equ:omega_as_equivalence}
  \begin{diagram}[labelstyle=\scriptstyle,midshaft,width=1.75cm]
    \GrTFInj(\A)  & \pile{ \rTo^{\pi} \\
                           \lTo_{\omega} } & \TailsInj(\A) \\
  \end{diagram}
\end{equation}
between the subcategory of torsion-free injective objects of $\Gr(\A)$
and the subcategory of all injective objects of $\Tails(A)$.

Let me next turn to derived categories.  The projection functor $\pi$
is exact so extends to a triangulated functor $\DGrA
\stackrel{\pi}{\longrightarrow} \DTailsA$ between derived categories.
Moreover, since $\A$ has finite global dimension, each object of the
category $\Gr(A)$ has a bounded resolution by injective objects.  The
same therefore holds for $\Tails(A)$, as one sees using $\omega$ and
$\pi$.  So right-derived functors can be defined on the unbounded
derived categories $\DGrA$ and $\DTailsA$ by
\cite[sec. 10.5]{Weibelbook}.

In particular, $\DTailsA \stackrel{\Romega}{\longrightarrow} \DGrA$
exists, and it is not hard to see that
\begin{equation}
\label{equ:pi_Romega}
  \begin{diagram}[labelstyle=\scriptstyle,midshaft,width=1.75cm]
    \DGrA & \pile{ \rTo^{\pi} \\
                   \lTo_{\Romega} } & \DTailsA \\
  \end{diagram}
\end{equation}
is an adjoint pair of functors.
\end{Remark}

\begin{Definition}
\label{dfn:k_perp}
Let
\[
  k^{\perp} = \{\, N \in \DGrA \,\mid\, \RHom_A(k,N) = 0 \,\}.
\]
\end{Definition}

\begin{Proposition}
\label{prp:Romega_as_equivalence}
The functors in equation \eqref{equ:pi_Romega} restrict to
a pair of quasi-inverse equivalences of triangulated categories
\begin{equation}
\label{equ:Romega_as_equivalence}
  \begin{diagram}[labelstyle=\scriptstyle,midshaft,width=1.75cm]
    k^{\perp} & \pile{ \rTo^{\pi} \\
                       \lTo_{\Romega} } & \DTailsA. \\
  \end{diagram}
\end{equation}
\end{Proposition}

\begin{proof}
First observe that diagram \eqref{equ:omega_as_equivalence} extends to
a pair of quasi-inverse equivalences
\begin{equation}
\label{equ:Romega_as_equivalence_proto}
  \begin{diagram}[labelstyle=\scriptstyle,midshaft,width=1.75cm]
    \sK(\GrTFInj\,\A)  & \pile{ \rTo^{\pi} \\
                                \lTo_{\omega} } & \sK(\TailsInj\,\A) \\
  \end{diagram}
\end{equation}
between the homotopy category of complexes of torsion free injective
objects of $\Gr(\A)$, and the homotopy category of complexes of
injective objects of $\Tails(\A)$.

Next, finite global dimension of $\A$ implies that $\DGrA$ is
equivalent to $\sK(\GrInj\,\A)$, the homotopy category of complexes of
injective objects of $\Gr(\A)$, and that under the equivalence, a
right-derived functor $\operatorname{R}\!F$ on $\DGrA$ corresponds to
the restriction of $F$ to $\sK(\GrInj\,\A)$.  See e.g.\
\cite[sec.\ 10.5]{Weibelbook}.  A similar remark applies to $\DTailsA$
and $\sK(\TailsInj\,\A)$.  Therefore, up to equivalence, diagram
\eqref{equ:pi_Romega} is
\begin{equation}
\label{equ:model1}
  \begin{diagram}[labelstyle=\scriptstyle,midshaft,width=1.75cm]
    \sK(\GrInj\,\A) & \pile{ \rTo^{\pi} \\
                             \lTo_{\omega} } & \sK(\TailsInj\,\A). \\
  \end{diagram}
\end{equation}

Translating diagram \eqref{equ:Romega_as_equivalence_proto} through
the equivalence between diagrams \eqref{equ:model1} and
\eqref{equ:pi_Romega} shows that diagram
\eqref{equ:Romega_as_equivalence_proto} gives an equivalence between
some subcategory of $\DGrA$, and the whole category $\DTailsA$.  To
finish the proof, I must show that the subcategory in question is
$k^{\perp}$.  That is, I must show that the subcategory
$\sK(\GrTFInj\,\A)$ of $\sK(\GrInj\,\A)$ corresponds to the
subcategory $k^{\perp}$ of $\DGrA$.

For this, note that by the above, the functor $\Hom_{\A}(k,-)$ on
$\sK(\GrInj\,\A)$ corresponds to the derived functor $\RHom_{\A}(k,-)$
on $\DGrA$, so I must show that $\sK(\GrTFInj\,\A)$ is the subcategory
of $\sK(\GrInj\,\A)$ annihilated by $\Hom_{\A}(k,-)$.

In fact, this is not quite true, but it is true and easy to see that
the subcategory of $\sK(\GrInj\,\A)$ annihilated by $\Hom_A(k,-)$
consists exactly of the complexes isomorphic to complexes in
$\sK(\GrTFInj\,\A)$, and this is enough. 
\end{proof}

\section{Koszul duality}
\label{sec:Koszul}

\begin{Remark}
\label{rmk:Koszul}
Let me recapitulate one of the versions of Koszul duality set up in
\cite{FloystadKoszul}.  According to \cite[thm.\
7.2.3']{FloystadKoszul}, there is a pair of quasi-inverse equivalences
of triangulated categories
\begin{equation}
\label{equ:pre_Koszul}
  \begin{diagram}[labelstyle=\scriptstyle,midshaft,width=1.75cm]
    \sK(\GrCoFree\,\Ad) & \pile{ \rTo^{\BGGF} \\
                                 \lTo_{\BGGG} } & \sK(\GrFree\,\A). \\
  \end{diagram}
\end{equation}
Here $\GrFree(\A)$ is the full subcategory of $\Gr(\A)$ consisting of
modules which have the form $\coprod_i \A(n_i)$, and $\GrCoFree(\Ad)$ is the
full subcategory of $\Gr(\Ad)$ consisting of modules which have the form
$\prod_j (\Ad)^{\prime}(m_j)$.  The categories
$\sK(\GrCoFree\,\Ad)$ and $\sK(\GrFree\,\A)$ are the corresponding
homotopy categories of complexes.

The functors $\BGGF$ and $\BGGG$ are constructed as follows in
\cite[sec.\ 3.2]{FloystadKoszul}.  Given $M$ in $\sK(\GrCoFree\,\Ad)$,
one constructs a double complex 
\[
  \begin{diagram}[objectstyle=\scriptstyle,labelstyle=\scriptscriptstyle,midshaft,width=0.8cm,height=0.8cm]
    & & \vdots & & \vdots & & \vdots & & \\
    & & \uTo & & \uTo & & \uTo & & \\
    \cdots & \rTo & \A(1) \otimes M_1^{-1}
           & \rTo &  \A(1) \otimes M_1^0 
           & \rTo &  \A(1) \otimes M_1^1 & \rTo & \cdots \\
    & & \uTo & & \uTo & & \uTo & & \\
    \cdots & \rTo & \A \otimes M_0^{-1}
           & \rTo &  \A \otimes M_0^0 
           & \rTo &  \A \otimes M_0^1 & \rTo & \cdots \\
    & & \uTo & & \uTo & & \uTo & & \\
    \cdots & \rTo & \A(-1) \otimes M_{-1}^{-1}
           & \rTo &  \A(-1) \otimes M_{-1}^0 
           & \rTo &  \A(-1) \otimes M_{-1}^1 & \rTo & \cdots \\
    & & \uTo & & \uTo & & \uTo & & \\
    & & \vdots & & \vdots & & \vdots & & \\
  \end{diagram}
\]
with certain differentials, and the total complex $\Tot^{\coprod}$,
defined using coproducts, is $\BGGF(M)$.  In the diagram, superscripts
indicate cohomological degree and subscripts indicate graded degree,
so for instance, the graded module in cohomological degree zero of the
complex $M$ is $M^0$.  Also, $\otimes$ indicates tensor product over
$k$.

And given $N$ in $\sK(\GrFree\,\A)$, one constructs a double complex
\[
  \begin{diagram}[objectstyle=\scriptstyle,labelstyle=\scriptscriptstyle,midshaft,width=0.8cm,height=0.8cm]
    & & \vdots & & \vdots & & \vdots & & \\
    & & \uTo & & \uTo & & \uTo & & \\
    \cdots & \rTo & (\Ad)^{\prime}(1) \otimes N_1^{-1}
           & \rTo &  (\Ad)^{\prime}(1) \otimes N_1^0 
           & \rTo &  (\Ad)^{\prime}(1) \otimes N_1^1 & \rTo & \cdots \\
    & & \uTo & & \uTo & & \uTo & & \\
    \cdots & \rTo & (\Ad)^{\prime} \otimes N_0^{-1}
           & \rTo &  (\Ad)^{\prime} \otimes N_0^0 
           & \rTo &  (\Ad)^{\prime} \otimes N_0^1 & \rTo & \cdots \\
    & & \uTo & & \uTo & & \uTo & & \\
    \cdots & \rTo & (\Ad)^{\prime}(-1) \otimes N_{-1}^{-1}
           & \rTo &  (\Ad)^{\prime}(-1) \otimes N_{-1}^0 
           & \rTo &  (\Ad)^{\prime}(-1) \otimes N_{-1}^1 & \rTo & \cdots \\
    & & \uTo & & \uTo & & \uTo & & \\
    & & \vdots & & \vdots & & \vdots & & \\
  \end{diagram}
\]
with certain differentials, and the total complex $\Tot^{\prod}$,
defined using products, is $\BGGG(N)$.

Finite global dimension of $\A$ implies that $\DGrA$ is equivalent to
$\sK(\GrFree\,\A)$ (see \cite[sec.\ 10.5]{Weibelbook}), so the
equivalences \eqref{equ:pre_Koszul} can also be read as
\begin{equation}
\label{equ:Koszul}
  \begin{diagram}[labelstyle=\scriptstyle,midshaft,width=1.75cm]
    \sK(\GrCoFree\,\Ad) & \pile{ \rTo^{\BGGF} \\
                                 \lTo_{\BGGG} } & \DGrA. \\
  \end{diagram}
\end{equation}
\end{Remark}

\begin{Remark}
The name Koszul duality is potentially confusing since ``duality''
might lead one to think of contravariant functors, while $\BGGF$ and
$\BGGG$ are in fact covariant.
\end{Remark}

For the following lemma, note that I use $\Sigma^i(-)$ for $i$'th
suspension, so if $M$ is a complex then $(\Sigma^i M)^{\ell} = M^{i +
\ell}$.

\begin{Lemma}
\label{lem:Koszul}
The functors $\BGGF$ and $\BGGG$ satisfy the following.
\begin{enumerate}

  \item  $\BGGF(M(i)) \cong \Sigma^i(\BGGF\!M)(-i)$.

  \item  $\BGGG(N(j)) \cong \Sigma^j(\BGGG\!N)(-j)$.

  \item  $\BGGF((\Ad)^{\prime})$ is isomorphic to the \Alm\ $k$,
         when both are viewed as objects of $\DGrA$.

\end{enumerate}
\end{Lemma}

\begin{proof}
(i) and (ii) can be seen by playing with the double complexes which
define $\BGGF$ and $\BGGG$.  (iii) holds by 
\cite[exam.\ 3.1.1]{FloystadKoszul}. 
\end{proof}

\begin{Remark}
\label{rmk:GrStab}
The injective stable category over a ring is defined as the
module category modulo the ideal of morphisms which factor through an
injective module.

The present paper uses the graded version of this, so the injective
stable category $\GrStab(\Ad)$ is defined as $\Gr(\Ad)$ modulo the
ideal of morphisms which factor through an injective object of
$\Gr(\Ad)$.

Since $\Ad$ is graded Frobenius by remark \ref{rmk:blanket}(v), the
category $\GrStab(\Ad)$ is triangulated.  For $M$ in $\GrStab(\Ad)$,
the suspension $\Sigma M$ is the first syzygy in an injective
resolution of $M$.  So $\Sigma M$ is the cokernel of an
injective pre-envelope, that is, an injective homomorphism $M
\longrightarrow I$ in $\Gr(\Ad)$ where $I$ is an injective object of
$\Gr(\Ad)$.  Note that any injective pre-envelope can be used;
changing the injective pre-envelope does not change the isomorphism
class of $\Sigma M$ in $\GrStab(\Ad)$.

The degree shifting functor $(-)(1)$ on $\Gr(\Ad)$ induces a
functor on $\GrStab(\Ad)$ which I will also denote $(-)(1)$.

It is not hard to prove that $\GrStab(\Ad)$ is equivalent to the full
subcategory of exact complexes in $\sK(\GrCoFree\,\Ad)$.  Under the
equivalence, a module $M$ corresponds to a complete cofree resolution
$C$ of $M$, that is, a complex $C$ in $\sK(\GrCoFree\,\Ad)$ which is
exact and has zeroth cycle module $\Z^0(C)$ isomorphic to $M$.  To
prove that this gives an equivalence, one uses that injective and
projective objects of $\Gr(\Ad)$ coincide because $\Ad$ is graded
Frobenius; in particular, the objects of $\GrCoFree(\Ad)$ are both
injective (since they are products of degree shifts of
$(\Ad)^{\prime}$) and projective, and hence it is possible to
construct complete cofree resolutions by splicing left-resolutions
with right-resolutions.

Under the equivalence between $\GrStab(\Ad)$ and the full subcategory
of exact complexes in $\sK(\GrCoFree\,\Ad)$, the suspension $\Sigma$
on $\GrStab(\Ad)$ corresponds to the ordinary suspension $\Sigma$ on
$\sK(\GrCoFree\,\Ad)$, given by moving complexes one step to the left
and switching signs of differentials.  Also, the functor $(-)(1)$ on
$\GrStab(\Ad)$ corresponds to the functor $(-)(1)$ on
$\sK(\GrCoFree\,\Ad)$ induced by degree shifting of \Adlm s.
\end{Remark}

\begin{Proposition}
\label{prp:Koszul}
The functors in equation \eqref{equ:Koszul} induce a pair of
quasi-inverse equivalences of triangulated categories
\[
  \begin{diagram}[labelstyle=\scriptstyle,midshaft,width=1.75cm]
    \GrStab(\Ad) & \pile{ \rTo \\
                          \lTo    } & k^{\perp}. \\
  \end{diagram}
\]
\end{Proposition}

\begin{proof}
Remark \ref{rmk:GrStab} identifies $\GrStab(\Ad)$ with the full
subcategory of exact complexes in $\sK(\GrCoFree\,\Ad)$, and
definition \ref{dfn:k_perp} defines $k^{\perp}$ as a full subcategory
of $\DGrA$.  To prove the proposition, I must show that these
subcategories are mapped to each other by the functors $\BGGF$ and
$\BGGG$ of equation \eqref{equ:Koszul}.

However, let $N$ be in $\DGrA$.  Then the $j$'th graded component of
the $i$'th cohomology module of the complex $\BGGG\!N$ is 
\begin{eqnarray*}
  \h^i(\BGGG\!N)_j 
    & \stackrel{\rm (a)}{\cong} 
    & \Hom_{\scrsK(\scrGr\,\Ad)}(\Ad,\Sigma^i(\BGGG\!N)(j)) \\
  & \stackrel{\rm (b)}{\cong} 
    & \Hom_{\scrsK(\scrGrCoFree\,\Ad)}((\Ad)^{\prime}(-d),\Sigma^i(\BGGG\!N)(j)) \\
  & \cong & (*),
\end{eqnarray*}
where (a) is classical and (b) holds because of $\Ad \cong
(\Ad)^{\prime}(-d)$, cf.\ remark \ref{rmk:blanket}(vi).  Adjointness
between $\BGGF$ and $\BGGG$ gives (c) in
\begin{eqnarray*}
  (*) & \stackrel{\rm (c)}{\cong}
    & \Hom_{\scrDGrA}(\BGGF((\Ad)^{\prime}),\Sigma^i N(j+d)) \\
  & \stackrel{\rm (d)}{\cong}
    & \Hom_{\scrDGrA}(k,\Sigma^i N(j+d)) \\
  & \cong
    & \h^i \RHom_{\A}(k,N)_{j+d},
\end{eqnarray*}
and (d) is by lemma \ref{lem:Koszul}(iii).

But now it is clear that $\BGGG\!N$ is exact if and only if $N$ is in
$k^{\perp}$, as desired.
\end{proof}

\section{The \BGGc}
\label{sec:BGG}

Composing the equivalences of categories from propositions
\ref{prp:Romega_as_equivalence} and \ref{prp:Koszul} gives the
following main theorem of the paper.

\begin{Theorem}
[The \BGGc]
\label{thm:BGG}
There are quasi-inverse equivalences of triangulated categories
\[
  \begin{diagram}[labelstyle=\scriptstyle,midshaft,width=1.75cm]
    \GrStab(\Ad) & \pile{ \rTo^{\BGGphi} \\
                          \lTo_{\BGGgamma}  } & \DTailsA. \\
  \end{diagram}
\]
\end{Theorem}

\begin{Example}
\label{exa:classical_BGG}
If $\A$ is the polynomial algebra $k[X_1,\ldots,X_d]$ then it is
classical that $\A$ satisfies the conditions of setup
\ref{set:blanket}, and the definition of $\Ad$ in setup
\ref{set:blanket} makes it easy to see that $\Ad$ is the exterior
algebra $\Ealg = \bigwedge(Y_1,\ldots,Y_d)$.  Also, $\Tails(\A)$ is
equivalent to the category $\QCoh(\BP^{d-1})$ of quasi-coherent
sheaves on $(d-1)$-dimensional projective space by Serre's theorem,
\cite[thm., p.\ 229]{ArtinZhang}.  So theorem \ref{thm:BGG} gives an
equivalence of categories
\[
    \GrStab(\Ealg) 
    \simeq \D(\QCoh\,\BP^{d-1}).
\]
This is the classical \BGGc, originally found in \cite[thm.\ 2]{BGG},
with the slight improvement of dealing with the stable category of all
modules and the unbounded derived category of quasi-coherent sheaves
rather than the finite subcategories dealt with in \cite[thm.\
2]{BGG}.
\end{Example}

\begin{Remark}
\label{rmk:BGG}
The functors $\BGGphi$ and $\BGGgamma$ from theorem \ref{thm:BGG} are
constructed by composing the functors from propositions
\ref{prp:Romega_as_equivalence} and \ref{prp:Koszul}.  Untangling this
gives the following concrete descriptions of $\BGGphi$ and $\BGGgamma$. 

To get $\BGGphi(M)$, first take a complete cofree resolution $C$ of
$M$ and look at $\BGGF(C)$, where $\BGGF$ is one of the functors from remark
\ref{rmk:Koszul}.  Then $\BGGF(C)$ is in $\sK(\GrFree\,\A)$, and may
also be viewed as being in $\DGrA$.  Now $\BGGphi(M) = \pi\!\BGGF(C)$,
where $\pi$ is one of the functors from proposition
\ref{prp:Romega_as_equivalence}.

To get $\BGGgamma(\cM)$, first look at $\Romega(\cM)$ where $\Romega$
is one of the functors from proposition
\ref{prp:Romega_as_equivalence}.  Then $\Romega(\cM)$ is in $\DGrA$,
and may also be viewed as being in $\sK(\GrFree\,\A)$.  So
$\BGGG(\Romega(\cM))$ is in $\sK(\GrCoFree\,\Ad)$, where $\BGGG$ is one
of the functors from remark \ref{rmk:Koszul}.  In fact,
$\BGGG(\Romega(\cM))$ is even in the full subcategory of exact
complexes in $\sK(\GrCoFree\,\Ad)$.  Now $\BGGgamma(\cM) =
\Z^0\!\BGGG(\Romega(\cM))$, where $\Z^0$ takes the zeroth cycle
module.
\end{Remark}

The next lemma follows immediately from lemma \ref{lem:Koszul}, parts
(i) and (ii).

\begin{Lemma}
\label{lem:BGG}
The functors $\BGGphi$ and $\BGGgamma$ satisfy the following.
\begin{enumerate}

  \item  $\BGGphi(M(i)) \cong \Sigma^i(\BGGphi\!M)(-i)$.

  \item  $\BGGgamma(\cM(j)) \cong \Sigma^j(\BGGgamma\!\cM)(-j)$.

\end{enumerate}
\end{Lemma}

For the following lemma, let $L$ be the minimal free resolution of the
graded \Adlm\ $k$.  Each $L^i$ is free and hence cofree because remark
\ref{rmk:blanket}(vi) says $(\Ad)^{\prime} \cong \Ad(d)$.  So $L$ is a
complex in $\sK(\GrCoFree\,\Ad)$, and I can apply the functor $\BGGF$
from remark \ref{rmk:Koszul} and get a complex $\BGGF(L)$ in $\DGrA$.

\begin{Lemma}
\label{lem:BGGFL_is_torsion}
The cohomology of $\BGGF(L)$ is torsion.
\end{Lemma}

\begin{proof}
By \cite[sec.\ 3.2]{FloystadKoszul}, the functor $\BGGF$ exists in an
alternative version, namely as a functor
\begin{equation}
\label{equ:F_on_Chain}
  \Chain(\Gr\,\Ad) \stackrel{\BGGF}{\longrightarrow} \Chain(\Gr\,\A),
\end{equation}  
where $\Chain$ denotes categories whose objects are complexes and
whose morphisms are chain maps.  This alternative version of $\BGGF$
induces the one from remark \ref{rmk:Koszul} because homotopy
categories of complexes can be obtained from categories $\Chain$ by
dividing away the ideals of null homotopic morphisms.

The version of $\BGGF$ in equation \eqref{equ:F_on_Chain} respects
small colimits by \cite[sec.\ 4.3]{FloystadKoszul}.  In
$\Chain(\Gr\,\Ad)$ the object $L$ is the colimit of the objects
\[
  L \langle j \rangle = 
  \cdots \longrightarrow 0 \longrightarrow L^{-j}
  \longrightarrow \cdots \longrightarrow L^0 
  \longrightarrow 0 \longrightarrow \cdots,
\]
so 
\begin{equation}
\label{equ:colim}
  \BGGF(L) 
  \cong \BGGF(\colim L \langle j \rangle)
  \cong \colim \BGGF(L \langle j \rangle).
\end{equation}

Now, $\Ad$ is Koszul by \cite[prop.\ 2.9.1]{BGS}, and $L$ is the
minimal free resolution of $k$ over $\Ad$, and so
\begin{equation}
\label{equ:Li}
  L^{-i} \cong \coprod \Ad(-i).
\end{equation}
This implies that $L^{-i}$ is concentrated in graded degrees $i,
\ldots, d+i$ because $\Ad$ is concentrated in graded degrees $0,
\ldots, d$ by remark \ref{rmk:blanket}(iv).  So the construction in
remark \ref{rmk:Koszul} says that $\BGGF(L \langle j
\rangle)$ is $\Tot^{\coprod}$ of a double complex whose non-zero part
can be sketched as
\begin{equation}
\label{equ:truncated_double_complex}
  \begin{diagram}[labelstyle=\scriptstyle,midshaft,width=0.8cm,height=0.8cm]
    \A(d+j) \otimes L_{d+j}^{-j} & & & & & & \\
    \uTo & & & & & & \\
    \A(d+j-1) \otimes L_{d+j-1}^{-j} & \rTo & \cdots & & & & \\
    \uTo & & & & & & \\
    \vdots & \rTo & \cdots & \rTo & \cdots & \rTo & \A(d) \otimes L_d^0 \\
    \uTo & & & & & & \uTo \\
    \A(j) \otimes L_j^{-j} & \rTo & \cdots & \rTo & \cdots & \rTo &  \vdots \\
    & & & & & & \uTo \\
    & & & & \cdots & \rTo & \A(1) \otimes L_1^0 \\
    & & & & & & \uTo \\
    & & & & & & \A \otimes L_0^0 \lefteqn{.} \\
  \end{diagram}
\end{equation}

Also, combining equation \eqref{equ:Li} with $\Ad \cong
(\Ad)^{\prime}(-d)$ which holds by remark \ref{rmk:blanket}(vi) gives
$L^{-i} \cong \coprod (\Ad)^{\prime}(-d-i)$.  So up to degree shift
and suspension, the $(-i)$'th column of
\eqref{equ:truncated_double_complex} is just a coproduct of copies of
the column obtained from $(\Ad)^{\prime}$.  This column has non-zero
part
\[
  \begin{diagram}[labelstyle=\scriptstyle,midshaft,width=0.8cm,height=0.8cm]
    \A \otimes (\Ad)^{\prime}_0 \\
    \uTo \\
    \vdots \\
    \uTo \\
    \A(-d) \otimes (\Ad)^{\prime}_{-d} \lefteqn{,} \\
  \end{diagram}
\]
and is a free resolution of the \Alm\ $k$, as follows from
\cite[exam.\ 3.1.1]{FloystadKoszul}.  So the columns of
\eqref{equ:truncated_double_complex} have cohomology only at the top
ends, and the cohomology in the $(-i)$'th column is $\coprod k(d+i)$.

Now consider the first spectral sequence of the double complex
\eqref{equ:truncated_double_complex} (see \cite[sec.\ 5.6]{Weibelbook}).
The previous part of the proof shows that the $E_2$-term of the
spectral sequence is non-zero only at the top ends of the columns of
\eqref{equ:truncated_double_complex}, where
\[
  E_2^{0d} \cong \coprod k(d), \; \ldots \;,
  E_2^{-j,d+j} \cong \coprod k(d+j).
\] 
Since the double complex is bounded in all directions, the spectral
sequence converges towards the cohomology of $\Tot^{\coprod}$.
Consequently, $\Tot^{\coprod}$ of the double complex has cohomology
only in cohomological degree $d$, and this cohomology sits in graded
degrees $-d, \ldots, -d-j$.

But this $\Tot^{\coprod}$ is $\BGGF(L \langle j \rangle)$.  So
equation \eqref{equ:colim} now shows that $\BGGF(L)$ has cohomology
only in cohomological degree $d$, and that this cohomology can be
non-zero only in graded degrees $-d,-d-1,\ldots$.  In particular, the
cohomology of $\BGGF(L)$ is torsion.
\end{proof}

Now consider the graded \Adlm\ $k$ viewed as an object of
$\GrStab(\Ad)$, and consider $\calO$, the ``structure sheaf'' in
$\Tails(\A)$ defined by $\calO = \pi(\A)$.  Then $\calO$ can also be
viewed as a complex in $\DTailsA$ concentrated in cohomological degree
zero, and the following result holds.

\begin{Theorem}
\label{thm:BGGphi_k}
The functor $\BGGphi$ satisfies $\BGGphi(k) \cong \calO$.
\end{Theorem}

\begin{proof}
To get $\BGGphi(k)$, I must take $\pi\!\BGGF(C)$, where $C$ is a
complete cofree resolution of the \Adlm\ $k$, while $\pi$ and $\BGGF$
are the functors from proposition \ref{prp:Romega_as_equivalence} and
remark \ref{rmk:Koszul} (cf.\ remark \ref{rmk:BGG}).

Let $L$ be a minimal free resolution of $k$, as in lemma
\ref{lem:BGGFL_is_torsion}.  Also, consider the functor $\BGGG$ from
remark \ref{rmk:Koszul}.  From \cite[exam.\ 3.1.1]{FloystadKoszul}
follows that $\BGGG(\A)$ is a cofree resolution of $k$.  So there are
canonical morphisms $L \longrightarrow k$ and $k \longrightarrow
\BGGG(\A)$ which compose to a morphism $L \longrightarrow \BGGG(A)$
whose mapping cone $C$ is easily seen to be a complete cofree
resolution of $k$.

The distinguished triangle $L \longrightarrow \BGGG(A)
\longrightarrow C \longrightarrow$ in $\sK(\GrCoFree\,\Ad)$
gives a distinguished triangle
\[
  \pi\!\BGGF(L) 
  \longrightarrow \pi\!\BGGF\!\BGGG(\A) 
  \longrightarrow \pi\!\BGGF(C) 
  \longrightarrow
\]
in $\DTailsA$.  Let me compute the three complexes here: The
cohomology of $\BGGF(L)$ is torsion by lemma
\ref{lem:BGGFL_is_torsion}, so $\pi\!\BGGF(L) \cong 0$.  And $\BGGF$
and $\BGGG$ are quasi-inverse equivalences, so $\BGGF\!\BGGG(\A)$ is
isomorphic to $\A$, so $\pi\!\BGGF\!\BGGG(\A) \cong \pi(\A) = \calO$.  

Finally, $\pi\!\BGGF(C)$ is $\BGGphi(k)$ as
mentioned above.  So the distinguished triangle reads
\[
  0 \longrightarrow \calO \longrightarrow \BGGphi(k) \longrightarrow,
\]
proving $\BGGphi(k) \cong \calO$.
\end{proof}

\section{Computations}
\label{sec:computations}

This section contains computations, some involving the \BGGc, which
will be used on periodic injective resolutions in section
\ref{sec:periodicity}. 

The following lemma is just a graded version of \cite[cor.\
2.5.4(ii)]{Bensonbook}. 

\begin{Lemma}
\label{lem:Benson_isomorphism}
Let $M$ be in $\GrStab(\Ad)$.  There are canonical isomorphisms
\[
  \Hom_{\scrGrStab(\Ad)}(k,\Sigma^i M)
  \longrightarrow \Ext_{\scrGr(\Ad)}^i(k,M) 
\]
for $i \geq 1$.
\end{Lemma}

\begin{Lemma}
\label{lem:Exts}
Let $M$ be in $\GrStab(\Ad)$ and consider $\cM = \BGGphi(M)$ in
$\DTailsA$.  Then
\[
  \Ext_{\scrGr(\Ad)}^i(k,M(-i+j)) 
  \cong \Ext_{\scrTails(\A)}^{j}(\calO,\cM(i-j))
\]
for $i \geq 1$ and each $j$.
\end{Lemma}

\begin{proof}
This is a simple computation:
\begin{eqnarray*}
  \Ext_{\scrGr(\Ad)}^i(k,M(-i+j))
  & \stackrel{\rm (a)}{\cong}
    & \Hom_{\scrGrStab(\Ad)}(k,\Sigma^i M(-i + j)) \\
  & \stackrel{\rm (b)}{\cong}
    & \Hom_{\scrDTailsA}(\BGGphi\!k,\BGGphi(\Sigma^i M(-i + j))) \\
  & \stackrel{\rm (c)}{\cong}
    & \Hom_{\scrDTailsA}(\calO,\Sigma^{j}\cM(i + j)) \\
  & = & \Ext_{\scrTails(\A)}^{j}(\calO,\cM(i + j)), \\
\end{eqnarray*}
where (a) is by lemma \ref{lem:Benson_isomorphism} and (b) is by the
\BGGc, theorem \ref{thm:BGG}, while (c) is by theorem
\ref{thm:BGGphi_k} and lemma \ref{lem:BGG}(i).  
\end{proof}

For the following lemma, observe that the finitely generated graded
modules form a full subcategory $\grstab(\Ad)$ of $\GrStab(\Ad)$, and
that the complexes which have bounded cohomology consisting of objects
from the category $\tails(\A)$ form a full subcategory $\DbTailsA$ of
$\DTailsA$.

\begin{Lemma}
\label{lem:finiteness}
The subcategories $\grstab(\Ad)$ and $\DbTailsA$ map to each other
under the \BGGc\
$\begin{diagram}[labelstyle=\scriptstyle,midshaft,width=1.75cm]
    \GrStab(\Ad) & \pile{ \rTo^{\BGGphi} \\
                          \lTo_{\BGGgamma}  } & \DTailsA. \\
 \end{diagram}$
\end{Lemma}

\begin{proof}
It is not hard to check that $\grstab(\Ad)$ consists of the objects of
$\GrStab(\Ad)$ which are finitely built from objects of the form
$k(i)$.

Similarly, $\DbTailsA$ consists of the objects of $\DTailsA$ which are
finitely built from objects of the form $\calO(j)$.

But under the \BGGc, $k(i)$ corresponds to $\Sigma^i
\calO(-i)$ by theorem \ref{thm:BGGphi_k} and lemma \ref{lem:BGG}(i),
so the present lemma follows.
\end{proof}

\begin{Lemma}
\label{lem:derived_Serre_vanishing}
Let $\cM$ be in $\DbTailsA$.  Then for $i \gg 0$ I have
\[
  \Ext_{\scrTails(\A)}^{j}(\calO,\cM(i-j))
  \cong \Hom_{\scrTails(\A)}(\calO,\h^{j}(\cM)(i-j))
\]
for each $j$, where $\h^j(\cM)$ is the $j$'th cohomology of $\cM$.
\end{Lemma}

\begin{proof}
The algebra $\A$ has global dimension $d$ by assumption, so
$\tails(\A)$ has cohomological dimension at most $d-1$ by
\cite[prop.\ 7.10(3)]{ArtinZhang}, so $\Ext_{\scrTails(\A)}^{\geq
d}(\calO,\cN) = 0$ holds for each $\cN$ in $\tails(\A)$.

Moreover, $\A$ is even \ASreg\ by assumption, so $\tails(\A)$
satisfies Serre vanishing by \cite[thms.\ 8.1(1) and 7.4]{ArtinZhang}.
That is, given $\cN$ in $\tails(\A)$ and given $p$ with $1 \leq p \leq
d-1$, I have $\Ext_{\scrTails(\A)}^p(\calO,\cN(r)) = 0$ for $r \gg 0
$.

Now, given $\cN$, I can kill all the finitely many
$\Ext_{\scrTails(\A)}^p(\calO,\cN(r))$ which might be non-zero by
choosing $r$ large enough.  That is, given $\cN$ in $\tails(\A)$, I
have
\begin{equation}
\label{equ:Serre_vanishing}
  \mbox{
    $\Ext_{\scrTails(\A)}^p(\calO,\cN(r)) = 0$ 
    for $r \gg 0$ and each $p \geq 1$.
       }
\end{equation}

There is a convergent spectral sequence
\[
  E_2^{pq} = \Ext_{\scrTails(\A)}^p(\calO,\h^q(\cM)(i - j))
  \Rightarrow \Ext_{\scrTails(\A)}^{p+q}(\calO,\cM(i - j))
\]
by \cite[5.7.9]{Weibelbook} (convergence because the cohomology
$\h(\cM)$ is bounded).  By assumption on $\cM$, the finitely many
non-zero $\h^q(\cM)$'s are in $\tails(\A)$.  So equation
\eqref{equ:Serre_vanishing} implies that for $i - j \gg 0$, the term
$E_2^{pq}$ is concentrated on the line $p = 0$.  So the spectral
sequence collapses and gives
\begin{equation}
\label{equ:derived_Serre_vanishing}
  \Hom_{\scrTails(\A)}(\calO,\h^q(\cM)(i - j)) \cong
  \Ext_{\scrTails(\A)}^q(\calO,\cM(i - j))
\end{equation}
for $i - j \gg 0$ and each $q$.

Now observe that the isomorphism \eqref{equ:derived_Serre_vanishing}
also holds for $q \gg 0$, simply because both sides are then zero.
For the left hand side, this holds because $\h(\cM)$ is
bounded.  For the right hand side, use that $\h(\cM)$ is bounded and
that $\tails(\A)$ has cohomological dimension at most $d-1$ by
\cite[prop.\ 7.10(3)]{ArtinZhang}.

So setting $q$ equal to $j$, the isomorphism
\eqref{equ:derived_Serre_vanishing} holds for $j \gg 0$, and for
other values of $j$ I can force $i - j \gg 0$ by picking $i \gg
0$, and then the isomorphism also holds.  That is,
\[
  \Hom_{\scrTails(\A)}(\calO,\h^{j}(\cM)(i - j)) \cong
  \Ext_{\scrTails(\A)}^{j}(\calO,\cM(i - j))
\]
for $i \gg 0$ and each $j$, proving the lemma.
\end{proof}

\section{Periodic injective resolutions}
\label{sec:periodicity}

This section shows how the \BGGc\ can be used to understand
periodicity of certain injective resolutions over exterior algebras as
a geometric phenomenon.

I also show an analogous \nc\ example with much more complicated
behaviour, due to the more intricate nature of \nc\ geometry.

\medskip
\noindent
{\bf The commutative case. }
The periodicity in question was discovered by Eisenbud in \cite[thm.\
2.2]{Eisenbud}.  Let $\Ealg$ be the exterior algebra
$\bigwedge(Y_1,\ldots,Y_d)$, and recall that $\gr(\Ealg)$ is the
category of finitely generated graded $\Ealg$-\lm s.

\begin{Theorem}
[Eisenbud]
\label{thm:Eisenbud}
Let $M$ in $\gr(\Ealg)$ be without injective direct summands, and suppose
that the Bass numbers
\[
  \mu^i(M) = \dim_k \Ext_{\Ealg}^i(k,M)
\]
are bounded for $i \geq 0$.  

Then the minimal injective resolution $I$ of $M$ is periodic with period
one in the following sense: Up to isomorphism, $I^i$ is $I^0(i)$ and
$\partial_I^{i+1}$ is $\partial_I^1(i)$.
\end{Theorem}

In other words, up to isomorphism and degree shift, all the $I^i$ and
all the $\partial_I^i$ are the same.  (In fact, Eisenbud worked with
minimal free resolutions, but using Matlis duality on his result gives
theorem \ref{thm:Eisenbud}.)

This phenomenon can be understood geometrically in a very simple way,
using the \BGGc: The module $M$ can be translated to a geometric
object on $\BP^{d-1}$, and since the Bass numbers of $M$ are bounded,
this object turns out to have zero dimensional support.  Therefore the
object is stable under twisting, that is, tensoring by
$\cO_{\BP^{d-1}}(1)$, and translating back, this gives that $M$ is
its own first syzygy, and periodicity of the minimal injective
resolution follows.

In more detail, let $\A$ be the polynomial algebra
$k[X_1,\ldots,X_d]$ so I am in the situation of example
\ref{exa:classical_BGG}.  In particular, $\Ad$ is the exterior algebra
$\Ealg = \bigwedge(Y_1,\ldots,Y_d)$, and $\Tails(\A)$ is equivalent to
$\QCoh(\BP^{d-1})$.  Let $M$ be in $\gr(\Ealg)$, and suppose that the
Bass numbers
\[
  \mu^i(M) = \dim_k \Ext_{\Ealg}^i(k,M)
\]
are bounded for $i \geq 0$.  

The \BGGc\ associates to $M$ the object $\cM = \BGGphi(M)$ of
$\D(\QCoh\,\BP^{d-1})$.  In fact, $\cM$ is even in
$\Db(\QCoh\,\BP^{d-1})$ by lemma \ref{lem:finiteness}, so only
finitely many of the cohomologies $\h^{\ell}(\cM)$ are non-zero, and
each $\h^{\ell}(\cM)$ is coherent.

For $i \geq 1$ I have
\begin{eqnarray*}
  \mu^i(M) & = & \dim_k \Ext_{\Ealg}^i(k,M) \\
  & \stackrel{\rm (a)}{=} & 
    \sum_{j} \dim_k \Ext_{\scrGr(\Ealg)}^i(k,M(-i+j)) \\
  & \stackrel{\rm (b)}{=} &
    \sum_{j} \dim_k 
    \Ext_{\scrQCoh(\BP^{d-1})}^{j}(\cO_{\BP^{d-1}},\cM(i-j)) \\
  & = & (*),
\end{eqnarray*}
where in (a), I am being clever by using the degree shift $-i+j$
instead of simply $j$, and where (b) is by lemma \ref{lem:Exts}.
And for $i \gg 0$ I have
\begin{equation}
\label{equ:mu_rewritten}
  (*) = \sum_{j} \dim_k 
        \Hom_{\scrQCoh(\BP^{d-1})}(\cO_{\BP^{d-1}},\h^{j}(\cM)(i-j))
\end{equation}
by lemma \ref{lem:derived_Serre_vanishing}.

So if $\mu^i(M)$ is bounded for $i \geq 0$, then each summand in
equation \eqref{equ:mu_rewritten} is bounded for $i \gg 0$.  That is,
\begin{equation}
\label{equ:Hilbert}
  \dim_k 
  \Hom_{\scrQCoh(\BP^{d-1})}(\cO_{\BP^{d-1}},\h^{j}(\cM)(\ell))
\end{equation} 
is bounded for each $j$ and $\ell \gg 0$.  However, this is now a
geometric statement: For $\ell \gg 0$, the polynomial growth rate of
the numbers in equation \eqref{equ:Hilbert} equals the dimension of
the support of $\h^{j}(\cM)$ on $\BP^{d-1}$, as follows from
\cite[thm.\ I.7.5]{Hartshorne}.  So it follows that each of the
finitely many non-zero $\h^{j}(\cM)$ has zero dimensional support; in
other words, the support is a finite collection of points.

Now suppose that the ground field $k$ is infinite.  Then it is
possible to pick a hyperplane $H$ in $\BP^{d-1}$ which is disjoint
from the support of each $\h^{j}(\cM)$.  To $H$ corresponds an
injection $\cO_{\BP^{d-1}}(1) \hookrightarrow \cO_{\BP^{d-1}}$ which
is an isomorphism away from $H$.  Tensoring over $\cO_{\BP^{d-1}}$
with $\cM$ gives a morphism $\cM \otimes \cO_{\BP^{d-1}}(1)
\stackrel{\mu}{\longrightarrow} \cM \otimes \cO_{\BP^{d-1}}$ 
and $\h^{j}(\mu)$ is $\h^{j}(\cM) \otimes
\cO_{\BP^{d-1}}(1) \longrightarrow \h^{j}(\cM) \otimes
\cO_{\BP^{d-1}}$.  However, this is an isomorphism for each $j$
because $\cO_{\BP^{d-1}}(1) \hookrightarrow \cO_{\BP^{d-1}}$ is an
isomorphism away from $H$ and hence an isomorphism on the support of
each $\h^{j}(\cM)$.  So $\mu$ is an isomorphism in
$\D(\QCoh\,\BP^{d-1})$, proving
\[
  \cM(1) \cong \cM.
\]
Under the \BGGc\ this gives $\BGGgamma(\cM(1)) \cong \BGGgamma(\cM)$,
and using $\BGGgamma(\cM) = \BGGgamma\!\BGGphi(M) \cong M$ and lemma
\ref{lem:BGG}(ii) this can be rearranged to
\begin{equation}
\label{equ:periodicity}
  \Sigma M \cong M(1)
\end{equation}
in $\GrStab(\Ealg)$.

In $\GrStab(\Ealg)$, the suspension $\Sigma M$ is computed as the
first syzygy of $M$ in an injective resolution, cf.\ remark
\ref{rmk:GrStab}.  So equation \eqref{equ:periodicity} shows that
in $\GrStab(\Ealg)$, this first syzygy is just $M$ itself, with a
degree shift of one.  It is possible to improve this with a few
remarks: First, if $M$ is without injective direct summands, then it
is not hard to show that the isomorphism \eqref{equ:periodicity} lifts
to hold in $\Gr(\Ealg)$, if $\Sigma M$ is obtained as the first syzygy
in a {\em minimal} injective resolution of $M$.  Secondly, the
assumption that $k$ is infinite can be dropped using
\cite[prop.\ 2.5.8]{EGAIV2}.

Iterating equation \eqref{equ:periodicity} now shows that in the
minimal injective resolution $I$ of $M$, the syzygy $\Sigma^i M$ is
simply $M(i)$.  Hence the module $I^i$ must be $I^0(i)$, and the
morphism $\partial_I^{i+1}$ must be $\partial_I^1(i)$.  So I have
recovered theorem \ref{thm:Eisenbud}.

\medskip
\noindent
{\bf The \nc\ case.}
In the above argument, the minimal injective resolution is periodic
with period one because points in $\BP^{d-1}$ are invariant under
twisting.  It is known that this invariance breaks down when one
passes to \nc\ analogues of $\BP^{d-1}$.

Here the twist can move points, and it is possible to have orbits of
length $n$, for any finite $n$, and orbits of infinite length.  So it
is obvious to expect that suitable \nc\ analogues of the above
argument might give examples of algebras $\Ad$, analogous to $\Ealg$,
and modules $M$ where $\mu^i(M)$ is bounded for $i \geq 0$, and yet
where the minimal injective resolution of $M$ is periodic with period
$n$, or aperiodic.  Indeed, this turns out to hold.

Note that the following example of this behaviour is already known
from \cite{Smith}.  But the present geometric view through the \BGGc\
is new.

\begin{Setup}
Suppose that the ground field $k$ is algebraically closed, and suppose
that $\elliptic$ is an elliptic curve over $k$ with a line bundle
$\bundle$ of degree $d$ and an automorphism $\auto$ given by
translation by a point of $\elliptic$.  To these data, \cite[sec.\
8]{Smith} associates a so-called Sklyanin algebra which satisfies the
standing assumptions from setup \ref{set:blanket}.  Let $\A$ be this
algebra.  
\end{Setup}

Observe that $\A$ is a \nc\ analogue of the polynomial algebra
$k[X_1,\ldots,X_d]$ and that hence, the Koszul dual $\Ad$ is a \nc\
analogue of the exterior algebra $\bigwedge(Y_1, \ldots, Y_d)$.

\begin{Setup}
The construction of $\A$ in \cite[sec.\ 8]{Smith} is so that
$\elliptic$ sits inside $\BP(\A_1^{\prime})$.  So each point $\point$
on $\elliptic$ is also a point in $\BP(\A_1^{\prime})$, that is, a one
dimensional subspace of $\A_1^{\prime}$.  This subspace has an
annihilator $\point^{\perp}$ in $\A_1$, and the graded \Alm\
$\pointmodule(\point) = \A/\A \point^{\perp}$ is a so-called point
module.  That is, it is cyclic, and each graded piece in non-negative
degrees is one dimensional.

Let me now write 
\[
  \cM(\point) = \pi(\pointmodule(\point)).
\]
This is an object of $\tails(\A)$, and I view it as a complex
concentrated in cohomological degree zero.  This complex is an object
of $\DTailsA$, so finally the \BGGc\ gives the object
\[
  M(\point) = \BGGgamma(\cM(\point))
\]
in $\GrStab(\Ad)$.  In fact, $\cM(\point)$ viewed as an object of
$\DTailsA$ is in the subcategory $\DbTailsA$, so lemma
\ref{lem:finiteness} says that $M(\point)$ is even in $\grstab(\Ad)$.

Observe that $M(\point)$ is only well-defined up to isomorphism in
$\GrStab(\Ad)$, so when looking at $M(\point)$ as a graded \Adlm, I
can drop any injective direct summands, and so assume that $M(p)$ is
without injective direct summands.
\end{Setup}

Let me start by pointing out the following property of the modules
$M(\point)$. 

\begin{Proposition}
\label{prp:Bass_bounded}
The Bass numbers $\mu^i(M(\point))$ are bounded for $i \geq 0$.
\end{Proposition}

\begin{proof}
For $i \gg 0$ I have
\begin{eqnarray*}
  \mu^i(M(\point)) & = & \dim_k \Ext_{\Ad}^i(k,M(\point)) \\
  & \stackrel{\rm (a)}{=} & 
    \sum_{j} \dim_k \Ext_{\scrGr(\Ad)}^i(k,M(\point)(-i+j)) \\
  & \stackrel{\rm (b)}{=} &
    \sum_{j} \dim_k 
    \Ext_{\scrTails(\A)}^{j}(\calO,\cM(\point)(i-j)) \\
  & \stackrel{\rm (c)}{=} & 
    \sum_{j} \dim_k 
    \Hom_{\scrTails(\A)}(\calO,\h^{j}(\cM(\point))(i-j)) \\
  & = & (*),
\end{eqnarray*}
where in (a), I am being clever by using the degree shift $-i+j$
instead of $j$, where (b) is by lemma \ref{lem:Exts}, and where
(c) is by lemma \ref{lem:derived_Serre_vanishing}.  However,
the complex $\cM(\point)$ is just the object $\cM(\point)$ placed in
cohomological degree zero, so
\[
  (*) = \dim_k \Hom_{\scrTails(\A)}(\calO,\cM(\point)(i)) = (**),
\]
and since $\cM(\point)$ is $\pi(\pointmodule(\point))$ and $i$ is
large, this is
\[
  (**) = \dim_k \pointmodule(\point)_i = 1
\]
by \cite[thm.\ 8.1(1) and prop.\ 3.13(2)]{ArtinZhang}, because the
algebra $\A$ is \ASreg.
\end{proof}

Now some computations with the $M(\point)$'s.

\begin{Lemma}
\label{lem:Mp_determines_p}
The module $M(\point)$ determines $\point$.
\end{Lemma}

\begin{proof}
It is certainly true that $M(\point)$ determines $\cM(\point) \cong
\BGGphi(M(\point))$.  In turn, $\cM(\point)$ 
determines the tail $\pointmodule(\point)_{\geq n}$ for $n \gg 0$,
because when viewing $\cM(\point)$ as an object of $\tails(\A)$, I
have
\begin{equation}
\label{equ:pointmodule_tail}
  \pointmodule(\point)_{\geq n} 
  \cong \omega \pi(\pointmodule(\point))_{\geq n}
  = \omega(\cM(\point))_{\geq n}
\end{equation}
for $n \gg 0$ by \cite[thm.\ 8.1(1) and prop.\ 3.13(2)]{ArtinZhang}.
But $\pointmodule(\point)_{\geq n}$ determines $\point$ by \cite[sec.\
8]{Smith}, so $M(\point)$ does too.
\end{proof}

\begin{Lemma}
\label{lem:tau}
The modules $M(\point)$ satisfy
\[
  \Sigma(M(\point)) \cong M(\auto^{2-d}\point)(1)
\]
in $\GrStab(\Ad)$.
\end{Lemma}

\begin{proof}
In \cite[exam.\ 9.5]{Smith} is proved
\[
  \pointmodule(\point)_{\geq 1}(1) \cong
  \pointmodule(\auto^{2-d}\point),
\]
and applying $\pi$ shows
\[
  \cM(\point)(1) \cong \cM(\auto^{2-d}\point)
\]
because $\pi$ only sees the tail of a module.  Applying $\BGGgamma$
and lemma \ref{lem:BGG}(ii), this can be rearranged to the lemma's
isomorphism
\[
  \Sigma(M(\point)) \cong M(\auto^{2-d}\point)(1).
\]
\end{proof}

\begin{Lemma}
\label{lem:shift}
If $\Sigma^i(M(\point)) \cong M(\pointtwo)(j)$ holds in $\GrStab(\Ad)$
for some points $\point$ and $\pointtwo$ on $\elliptic$, then $i = j$.
\end{Lemma}

\begin{proof}
The lemma's isomorphism implies $\BGGphi(\Sigma^i(M(\point))) \cong
\BGGphi(M(\pointtwo)(j))$, and using lemma \ref{lem:BGG}(i) and
$\BGGphi(M(\point)) = \cM(\point)$, this becomes
$\Sigma^i(\cM(\point)) \cong \Sigma^j(\cM(\pointtwo))(-j)$.  Since the
cohomologies of $\cM(\point)$ and $\cM(\pointtwo)$ are concentrated in
cohomological degree zero, this is only possible with $i = j$.
\end{proof}

Finally, these lemmas can be used as follows.  If there is to be
periodicity in the sense
\begin{equation}
\label{equ:periodicity2}
  \Sigma^i(M(\point)) \cong M(\point)(j) 
\end{equation}
in $\GrStab(\Ad)$ for some $i$ and $j$, then $i = j$ by lemma
\ref{lem:shift}.  On the other hand, $\Sigma^i(M(\point)) \cong
M(\auto^{(2-d)i}\point)(i)$ holds by lemma \ref{lem:tau}.  Substituting
into equation \eqref{equ:periodicity2} gives
$M(\auto^{(2-d)i}\point)(i) \cong M(\point)(i)$, hence
$M(\auto^{(2-d)i}\point) \cong M(\point)$, and as $M(\point)$
determines $\point$ by lemma \ref{lem:Mp_determines_p}, this implies
$\auto^{(2-d)i}(\point) = \point$. 

Conversely, $\auto^{(2-d)i}(\point) = \point$ gives
\[
  \Sigma^i(M(\point)) \cong M(\auto^{(2-d)i}\point)(i)
  \cong M(\point)(i)
\]
in $\GrStab(\Ad)$.

Summing up, if $d$, $\auto$ and $\point$ are so that
$\auto^{(2-d)i}(\point) \not= \point$ for $i = 1, \ldots, n-1$ but
$\auto^{(2-d)n}(\point) = \point$, then in $\GrStab(\Ad)$ the suspension
$\Sigma^i(M(\point))$ is not a degree shift of $M(\point)$ for $i = 1,
\ldots, n-1$, but $\Sigma^n(M(\point))$ is $M(\point)(n)$.

And if $d$, $\auto$ and $\point$ are so that $\auto^{(2-d)i}(\point)
\not= \point$ for $i \geq 1$, then in $\GrStab(\Ad)$ the suspension
$\Sigma^i(M(\point))$ is not a degree shift of $M(\point)$ for $i
\geq 1$.

Using that $M(\point)$ contains no injective direct summands, this
easily lifts to give the same result in $\Gr(\Ad)$ for syzygies in
{\em minimal} injective resolutions.  So I get the following example
which shows the promised contrast to theorem \ref{thm:Eisenbud} with
respect to periodicity of minimal injective resolutions.

\begin{Example}
\label{exa:periodicity}
(1)  Let $d$, $\auto$ and $\point$ be so that $\auto^{(2-d)i}(\point)
\not= \point$ for $i = 1, \ldots, n-1$ but $\auto^{(2-d)n}(\point) =
\point$. 

Then the minimal injective resolution $I$ of $M(\point)$ is periodic
with period $n$, in the sense that in the resolution, the $i$'th
syzygy $\Sigma^i(M(\point))$ is not isomorphic to a degree shift of
$M(\point)$ for $i = 1, \ldots, n-1$, but the $n$'th syzygy
$\Sigma^n(M(\point))$ is isomorphic to $M(\point)(n)$.

Hence up to isomorphism, $I^n$ is $I^0(n)$ and $\partial_I^{n+1}$ is
$\partial_I^1(n)$, while the same is not true with any smaller value
of $n$.

\smallskip
\noindent
(2)  Let $d$, $\auto$ and $\point$ be so that $\auto^{(2-d)i}(\point)
\not= \point$ for $i \geq 1$.

Then the minimal injective resolution $I$ of $M(\point)$ is aperiodic,
in the sense that in the resolution, no syzygy $\Sigma^i(M(\point))$
is a degree shift of $M(\point)$ for $i \geq 1$.
\end{Example}

\medskip
\noindent
{\bf Acknowledgement. } The diagrams were typeset with Paul Taylor's
{\tt diagrams.tex}.

\end{document}